\documentstyle[12pt]{article}

  \newcommand{\const}{\rm const}

  \newcommand{\Law}{\rm Law}

  \newcommand{\Poisson}{\rm  Poisson}
  \newcommand{\Ent}{\rm Ent}

\begin{document}

   \begin{center}

 \ {\bf Non-asymptotic estimation for Bell function, }\\

 \vspace{4mm}

 \  {\bf with probabilistic applications} \\

 \vspace{4mm}

\vspace{7mm}

  {\bf   Ostrovsky E., Sirota L.}\\

\vspace{4mm}

 Israel,  Bar-Ilan University, department  of Mathematic and Statistics, 59200, \\

\vspace{4mm}

e-mails:\ eugostrovsky@list.ru \\
sirota3@bezeqint.net \\

\vspace{4mm}

  {\bf Abstract} \\

\vspace{4mm}

 \end{center}

 \ \  We  deduce the non-asymptotical bilateral estimates for moment inequalities for sums of non-negative independent random variables,
based on the correspondent estimates for the so-called Bell functions and the Poisson distribution.\par

 \vspace{4mm}

 \par

\vspace{4mm}


\vspace{4mm}

 \section{ Definitions. Notations. Previous results. Statement of problem.}

\vspace{3mm}

 {\bf Definition 1.1.} The famous Bell {\it numbers} $ \  B(p), \ p \ge 0 \ $ are defined by means of the series

$$
B(p) \stackrel{def}{=} e^{-1} \ \sum_{k=0}^{\infty} \frac{k^p}{k!}. \eqno(1.0)
$$
 \ These numbers was introduces originally for the integer positive values $ \ p \ $ by E.T.Bell [1], [2]; see also Dobinski [6]. \par
 \ They plays a very important role in the combinatorics [21], theory of functions, asymptotical analysis [4] and
especially in the theory probability, in the theory of summation of independent random variables [3], [7], [8], [9]-[11], [12]-[14], [16], [17], [19], [24]-[26] etc. \par

 \  More generally, define the so-called Bell's {\it function} of two variables

$$
B(p,\beta) \stackrel{def}{=} e^{-\beta} \sum_{k=0}^{\infty} \frac{k^p \ \beta^k}{k!}, \ p \ge 2, \ \beta > 0, \eqno(1.1)
$$
so that $ \ B(p) = B(p,1).  \ $  On the other words,  $ \ B(p,\beta) \ $ are the Bell function depending on additional parameter $ \  \beta \in (0, \infty).  \  $ \par

 \ Let the random variable (r.v.)  $ \ \tau = \tau[\beta], \ $ defined on certain probability space $ \ (\Omega, F,{\bf P}) \ $ with expectation $ \ {\bf E}, \ $
 has a Poisson distribution with parameter $ \  \beta, \ \beta > 0; \   $ write $ \ \Law(\tau) = \Law \tau[\beta] = \Poisson(\beta):   \  $

$$
{\bf P}(\tau = k) = e^{-\beta} \frac{\beta^k}{k!}, \ k = 0,1,2,\ldots,
$$

 \ It is worth to note that

$$
B(p,\beta) = {\bf E} (\tau[\beta])^p, \ p \ge 0.\eqno(1.2)
$$

\ In detail, let $ \ \eta_j, \ j = 1,2,\ldots \ $ be a sequence of non-negative independent random (r.v.); the case of centered or moreover symmetrical distributed r.v.
was considered in many works, see e.g. [5], [9]-[11],  [15], [17], [19], [23], [24], [25], [26], and so one. \par
 \  The following inequality holds true

$$
{\bf E} \left( \sum_{j=1}^n \eta_j   \right)^p \le B(p) \ \max \left\{  \ \sum_{j=1}^n {\bf E} \eta_j^p, \ \left( \ \sum_{j=1}^n {\bf E} \eta_j \ \right)^p \ \right\}, \ p \ge 2, \eqno(1.3)
$$
where the "constant" $ \ B(p) \ $ in (1.2) is the best possible, see [3], [24]. \par

 \ One of the interest applications of these estimates in statistics, more precisely,  in the theory of $ \ U \ $ statistics  may be found in the article [8]. \par

 \ Another application. Let $ \  n = 1,2,3,\ldots; \  a, b = \const > 0; \ p \ge 2, \ \mu = \mu(a,b;p) := a^{p/(p-1)} \ b^{1/(1-p)} . \ $
 Define the following class of the sequences of an independent non-negative random variables

$$
Z(a,b) \stackrel{def}{=} \left\{ \eta_j, \ \eta_j \ge 0, \ \sum_{j=1}^n {\bf E}\eta_j = a; \ \sum_{j=1}^n {\bf E}\eta_j^p = b \ \right\}. \eqno(1.4)
$$

 \ G.Schechtman in [24] proved that

$$
\sup_{ \ n = 1,2,\ldots; \  \{\eta_j \} \in Z(a,b) \  } {\bf  E}  \left( \ \sum_{j=1}^n \eta_j \  \right)^p = \left(  \frac{b}{a} \right)^{p/(p-1)} \ B(p, \mu(a,b;p)). \eqno(1.5)
$$

 \ N. G. de Bruijn  in the book [4] proved the following  {\it asymptotical} logarithmical expression for $ \ B(p) \ $ as $ \ p \to \infty: $

$$
\frac{\ln B(p)}{p} = \ln p - \ln \ln p - 1 + \frac{\ln \ln p}{\ln p} +
$$

$$
\frac{1}{\ln p} + \frac{1}{2} \left( \ \frac{\ln \ln p}{\ln p} \ \right)^2 + 0 \left( \ \frac{\ln \ln p}{\ln^2 p} \ \right), \eqno(1.6)
$$
at last for the integer values $ \ p. \ $\par

\vspace{4mm}

 \ {\bf Our aim in this short report is obtaining the bilateral non-asymptotical estimates for introduced before Bell functions, with "constructive" values of constants.  }\par

\vspace{4mm}

 \ We refine ones in [3], [4] etc. \par

 \ Denote as ordinary for arbitrary numerical valued r.v. $ \eta \ $ its Lebesgue-Riesz $ \ L(p) \ $ norm by $ \ |\eta|_p: $

 $$
 |\eta|_p \stackrel{def}{=} \left[ {\bf E} |\eta|^p \right]^{1/p}, \ p \in [1, \infty).  \eqno(1.7)
 $$

 \vspace{4mm}

 \ There are also many works  about  problems raised here, see, e.g. [3], [5], [7], [9]-[11], [12], [15], [17], [19], [23], [24], [25], [26]  and other
 articles mentioned in the references. \par

\vspace{4mm}

 \section{ General approach. Upper and lower estimates.}

\vspace{4mm}

{\bf A. An upper estimate.} \par

\vspace{4mm}

 \ Let again  $ \ \Law(\tau) = \Law \tau[\beta] = \Poisson(\beta).   \  $  In order to estimate the moments of this variable, we will apply the theory of
the so-called Grand Lebesgue Spaces (GLS) \ [17]. \par
  \ Indeed, let us calculate the moment generating function (MGF) for this variable; it has a form
%&
$$
{\bf E} e^{\lambda \tau} = e^{-\beta} \sum_{k=0}^{\infty} \frac{\beta^k \ \exp(\lambda \ k)}{k!} = \exp \left( \beta \left(e^{\lambda} - 1 \right)   \right). \eqno(2.1)
$$

 \ We will use an elementary inequality

$$
x^p \le \left(  \frac{p}{e \ \lambda}  \right)^p \ e^{\lambda x}, \ \lambda, x > 0;
$$
 and have

$$
{\bf E} \xi^p \le \left(  \frac{p}{e \ \lambda}  \right)^p \exp \left[ \ \beta  \  \left(e^{\lambda} - 1 \ \right) \ \right],\eqno(2.2)
$$
or equally

$$
|\xi|_p  \le \frac{p}{e \ \lambda} \ \exp \left\{ \beta \ (e^{\lambda} - 1)/p  \right\}. \eqno(2.2a)
$$

 \ Let us introduce the following auxiliary function

$$
g_{\beta}(p) \stackrel{def}{=} \frac{p}{e} \cdot
 \inf_{\lambda > 0} \left[\lambda^{-1} \ \exp \left( \ \beta \left(e^{\lambda} - 1  \right) \ \right) \ \right]^{1/p}, \ \beta, p  > 0. \eqno(2.3)
$$

 \ As long as the relations (2.2), (2.2a) are true for arbitrary positive value $ \ \lambda, \ $ one can select as the value $ \ \lambda \ $ its  optimal value. We
obtained  actually a following upper estimate for the Bell function.  \par

\vspace{4mm}

{\bf Proposition 2.1.}

$$
B^{1/p}(p, \beta) \le  g_{\beta}(p), \ p,\beta > 0. \eqno(2.4)
$$

\vspace{6mm}

{\bf A. A lower estimate.} \par

\vspace{4mm}

 \ We have for the r.v. $ \ \tau = \tau_{\beta}: \ \Law \left(\tau_{\beta} \right)  = \Poisson(\beta) $

$$
{\bf E} \tau^p = e^{-\beta} \sum_{k=0}^{\infty} \frac{k^p \ \beta^k}{k!}, \ p, \beta > 0.
$$

 \ Let us introduce  the following function

$$
h_0(p, \beta) \stackrel{def}{=} \sup_{k=1,2,\ldots}  e^{-\beta}  \left\{ \ \frac{k^p \ \beta^k}{k!} \ \right\}; \eqno(2.5)
$$
therefore

$$
B(p,\beta) \ge h_0(p, \beta), \ p, \ \beta > 0. \eqno(2.6)
$$

\vspace{4mm}

 \ The last estimate may be simplified as follows. We will apply the following version of the famous Stirling's formula [22]

$$
 k! \le \zeta(k), \ \ k = 1,2,\ldots,
$$
where

$$
\zeta(k) \stackrel{def}{=} \sqrt{2 \pi k} \ \left( \frac{k}{e}  \right)^k \ e^{1/(12k)}, \ k = 1,2,\ldots \eqno(2.7)
$$

 \ It is worth to note that the function $ \ \zeta = \zeta(k) \ $ may be extended as the function of   the {\it real } variable $ \ x \in [1, \infty): \ \zeta = \zeta(x) \ $
by means of the formula (2.7). \par

 \ Define a new function

$$
h(p, \beta) \stackrel{def}{=} \sup_{x \in(1, \infty)} \
 \left\{ e^{ 1/(6 p x )} \cdot \left[ \ \frac{e^{x - \beta} \ x^{p - x - 1/2}}{ \sqrt{2 \ \pi} \ x^x}  \ \right]^{1/p} \right\}. \eqno(2.8)
$$

 \ We obtained really the following lower estimate for the Bell's function.\par

\vspace{4mm}

 \ {\bf  Proposition 2.2.  }

$$
B^{1/p} (p, \beta) \ge h_0(p, \ \beta), \  B^{1/p} (p, \beta) \ge h(p, \ \beta), \ p, \ \beta > 0. \eqno(2.9)
$$

\vspace{4mm}

 \ We concretize further in the next sections the choosing of the values $ \ k_0, \ x_0, \ \lambda \ $ in order to simplify
the estimates (2.4) \ and (2.9).\par

\vspace{4mm}

 \section{ Main result. Simplification of the upper estimate. The case of one variable.}

\vspace{4mm}

 \ Suppose here that $ \ \beta = \const > 0, \ p \ge 2 \beta, \ p \ge 1. \ $ One can choose in (2.3) the  (asymptotically as $ \ p \to \infty\ $  optimal) value

$$
\lambda := \lambda_0 \stackrel{def}{=} \ln (p/\beta) - \ln \ln (p/\beta). \eqno(3.1)
$$

 \ We deduce from the proposition 2.1 after substituting and some cumbersome calculations \par

\vspace{4mm}

{\bf Proposition 3.1.} We assert under our conditions $ \ \beta = \const > 0, \ p \ge 2 \beta, \\ \ p \ge 1 \ $

$$
B^{1/p}(p,\beta) \le \frac{p/e}{\ln (p/\beta) - \ln \ln (p/\beta)} \cdot \exp \left\{ \frac{1}{\ln (p/\beta)}  - \frac{1}{p/\beta}  \right\}. \eqno(3.2)
$$

 \ Notice that the expression in the right-hand of (3.2) is in strict accordance with the strict asymptotic for Bell's number (1.6) obtained by
N.G.de Bruijn in the book [4], still in the case when $ \ \beta = 1. \ $ \par

\vspace{4mm}

 \ For example,

$$
B^{1/p}(p) \le \frac{p/e}{\ln p - \ln \ln p} \cdot \exp \left\{ 1/\ln p - 1/p   \right\}, \ p \ge 2. \eqno(3.3)
$$

\vspace{4mm}

\ {\bf Remark 3.1.} The obtained estimation (3.2) may be rewritten as follows

$$
B^{1/p}(p) \le \frac{p}{e \ \ln (p/\beta)} \cdot \left[ 1 + C_1(\beta) \cdot \ \frac{\ln \ln (p/\beta)}{\ln (p/\beta)} \  \right], \eqno(3.4)
$$
where
$  \ C_1(\beta) = \const \in (0,\infty), $ and we recall that $ \ p \ge 2\beta;\ $
with "constructive" and absolute value of the "constant" $ \ C_1(\beta). \ $ \par

\vspace{4mm}

 \section{ Main result. Simplification of the lower estimate. The case of one variable.}

\vspace{4mm}

 \ Let us return to the relations (2.5), (2.6). One can choose the following value as a capacity of the number  $ \ k; \  k := k_0(p,\beta) \ $

 $$
 k_0(p, \beta) := \frac{p}{\ln (pe/\beta)}. \eqno(4.1)
 $$

 \ More precisely,

$$
 k_0(p, \beta) := \Ent \left[ \ \frac{p}{\ln (pe/\beta)} \ \right] + 1, \eqno(4.1a)
$$
where $ \ \Ent[z] \ $ denotes the integer part of the (positive) number $ \ z. \ $\par

 \ We get again after cumbersome calculations

 $$
 \ B^{1/p}(p,\beta)  \ge \
 $$

$$
 \beta^{1/\ln(pe/\beta)} \cdot \frac{p}{\ln (pe/\beta)} \cdot \left\{ \exp \left[ \ \frac{\ln p - \ln(pe)/\beta}{\ln(pe)/\beta} \ \right] \  \right\}^{-1},
$$

$$
 \ p,\beta > 0, p/\beta \ge 2. \eqno(4.2)
$$

 \ After simplifications:

$$
B^{1/p}(p, \beta) \ge \frac{p}{e \ \ln (p/\beta)} \cdot \left[ 1 - C_2(\beta) \cdot \ \frac{\ln \ln (p/\beta)}{\ln (p/\beta)} \  \right], \eqno(4.3)
$$
where
$  \ C_2(\beta) = \const \in (0,\infty), $ and we recall that $ \ p \ge 2\beta;\ $
with "constructive" and absolute value of the "constant" $ \ C_2(\beta) \in (0,\infty). \ $ \par

 \ Notice that the upper and lower bounds almost coincides and almost coincides with the asymptotic expression for Bell's numbers. \par

 \ As a slight corollary: under condition (5.2)

$$
\left| \ \frac{ \  B^{1/p}(p, \beta) - \frac{p}{e \ \ln (p/\beta)} \ }{ \frac{p}{e \ \ln (p/\beta)}\ } \ \right|
\le C_0(\beta) \frac{\ln \ln (p/\beta)}{\ln (p/\beta)}. \eqno(4.4)
$$

\vspace{4mm}

 \section{ Main result. The case of two variables. Upper bounds. }

\vspace{4mm}

  \ We give first of all a rough estimate for the Bell function. Namely, let as before the r.v. \\
 $ \ \tau = \tau[\beta] \ $ has the Poisson distribution with a parameter $  \ \beta; \beta > 0: \  $

 $$
  \ \Law(\tau) = \Poisson(\beta).  \eqno(5.0)
 $$

 \vspace{4mm}

 \ {\bf Proposition 5.1.}

$$
|\tau[\beta]|_p =
B^{1/p}(p,\beta) \le \beta \ \frac{p}{e \ \ln p} \ \left( \ 1 + C_3 \frac{\ln \ln p}{\ln p}\ \right), \ p \ge 2. \eqno(5.1)
$$
with some absolute constant $ C_3. \ $ \par
 \ The last estimate in (5.1) is essentially non-improvable at last in the case when $ \ \beta = 1. $ \par

\vspace{4mm}

 {\bf Proof.} We can and will suppose without loss of generality that the number $ \ \beta \ $ is integer. Introduce on an appropriate (sufficiently rich)
  probability space the sequence $ \ \{ \theta(i) \}, \ i = 1,2,\ldots \ $  of independent  standard Poisson distributed random variables:

$$
\Law(\theta(i)) = \Poisson(1), \ \Leftrightarrow {\bf P}(\theta(i) = k) = e^{-1} /k!, \ k = 0,1,2,\ldots.
$$

 \ The distribution of the  sum $ \ \sum_{i=1}^{\beta} \theta(i) \ $ coincides with $ \ \tau[\beta]; \ $ one can assume

$$
\tau[\beta] = \sum_{i=1}^{\beta} \theta(i).
$$
 \ One can apply a triangle inequality for the Lebesgue - Riesz norm $ \  L_p: \ $

$$
B^{1/p}(p,\beta) = |\tau[\beta]|_p  \le \sum_{i=1}^{\beta} |\theta(i)|_p = \beta \ |\theta(1)|_p.
$$
  \ It remains to use the proposition (3.1.) \par

\vspace{4mm}

 \ We suppose hereafter that both the variables $ \ p \ $ and $ \ \beta \ $ are independent but such that

$$
p \ge 1, \ \beta > 0, \ p/\beta  \le 2. \eqno(5.2)
$$

 \ It is this case namely that takes place in the work of G.Schechtman [24],  see (1.4). Indeed, suppose for simplicity therein that
the non- negative random variables $ \ \{\eta_j\}, \ i =1,2,\ldots; \ \eta := \eta_1, \ $ are independent and identical distributed (i,; i.d.) and such that for some $ \ p > 1 \ $

$$
m_1 := {\bf E} \eta < \infty; \ \ m_p := {\bf E} \eta^p < \infty.
$$
 \ Then in (1.4)

$$
a = n \ m_1,\ b = n \ m_p,
$$
so that $ \  \beta \asymp n, \ n \to \infty.  \  $\par

 \ We return to the estimate (2.2a)

$$
|\tau[\beta]|_p  \le \frac{p}{e \ \lambda} \ \exp \left\{ \beta \ (e^{\lambda} - 1)/p  \right\}. \eqno(5.3)
$$

 \  But one can now choose in (5.3) under our conditions the value  $ \ \lambda := p/\beta. \ $  \par

\vspace{4mm}

{\bf Proposition 5.2.} \ We get after simple calculations under formulated before in this section conditions  (5.2)

$$
B^{1/p}(p, \beta) \le  K_+ \cdot \beta,
$$
where

$$
 K_+  := \exp \left[ (e^2 - 3)/2   \right]  \approx 8.9758...  \eqno(5.4)
$$

\vspace{4mm}

 \section{ Main result. The case of two variables. Lower bounds. }

\vspace{4mm}

 \ We retain the conditions (5.2); and we will use also as above the following lower estimate

$$
B(p, \beta) \ge e^{-\beta} \cdot \frac{\beta^p \ p^p}{\sqrt{2 \ \pi}} \ k^{ -k + 1/2} \ e^{-k + 1/(12 k)}. \eqno(6.1)
$$

 \ One can choose in (6.1) the value $ \ k =  k_0 = k_0(p)  := p, \ $ if $ \ p \ $  is integer, and $ \ k_0 := \Ent(p) \ $ otherwise. We deduce
 after some calculations: \par

 \vspace{4mm}

 {\bf Proposition 6.1.} We assert  under conditions (5.2)

\vspace{4mm}

$$
B^{1/p}(p,\beta) \ge K_- \cdot \beta,
$$
where

$$
 K_- = (2 \pi)^{-1/2} \ \exp \left[  -1/(2e) \ + 1/3   \right] \ \approx 0.6538... \ . \eqno(6.2)
$$

\vspace{4mm}

 \section{ Concluding remarks. }

\vspace{4mm}

 \ {\bf A.} It is interest by our opinion to compute the exact value of the constants $ \ K_{\pm} \ $ in the estimates (5.4) and (6.2),
 as well as to find  its "limit" behavior.\par

\vspace{4mm}

\ {\bf B.} It is interest also by our opinion to generalize the approximation of the Bell function $ B(p) \ $  through  the so - called Lambert function $ \  W(\cdot), \ $
where by definition

$$
W(p) \ e^{W(p)} = p, \ p = 0,1,2,\ldots:
$$

$$
B(p) \sim \frac{1}{\sqrt{p}} \cdot \left( \ \frac{p}{W(p)} \ \right) \cdot \exp \left( \ \frac{p}{W(p)} - p - 1 \ \right)
$$
onto the more general function $ \ B(p,\beta). \ $ \par

\vspace{4mm}

 \begin{center}

 \vspace{6mm}

 {\bf References.}

 \vspace{4mm}

\end{center}

{\bf 1. Bell, E. T.} (1934). {\it Exponential polynomials.}  Annals of Mathematics. 35: 258-277. doi:10.2307/1968431. JSTOR 1968431. \\

\vspace{3mm}

{\bf 2. Bell, E. T.} (1938). {\it The iterated exponential integers. } Annals of Mathematics. 39: 539-557. doi:10.2307/1968633. JSTOR 1968633. \\

\vspace{3mm}

{\bf 3. Daniel Berend, Tamir Tassa.  } {\it Improved bounds on Bell numbers and on moments of sums of random variables.}
Probability and Mathematical Statistics, Vol. 30, Fasc. 2 (2010), pp. 185-205. \\

\vspace{3mm}

{\bf 4. N.G.de Bruijn.} {\it Asymptotic Methods in Analysis,}  Dover, New York, NY, 1958.\\

\vspace{3mm}

{\bf 5. Dharmadhikari S., Jogdeo K.} {\it Bounds  on  the  Moment of certain random Variables.}
Ann. Math. Statist., 1969, V. 40, B.4 pp. 1506-1518 \\

\vspace{3mm}

{\bf  6. G.Dobinski.} {\it  Summierung der reihe } $ \ e^{-1} \sum_n n^m/n! \ $ f\"ur $ \ m = 1,2,3,4,5,6,\ldots \ $.
Grunert  Archiv,  (Arch. Math. Phys.) {\bf 61,} (1877), 333-336.\\

\vspace{3mm}

{\bf 7. T. Figiel, P. Hitczenko, W. B. Johnson, G. Schechtman and J. Zinn.} {\it Extremal
properties of Rademacher functions with applications to the Khinchine and Rosenthal inequalities,}
Trans. American. Math. Society. 349, (1997), pp. 997-1027.

\vspace{3mm}

{ \bf 8. Evarist Gine, Rafal Latala, Joel Zinn.}  {\it Exponential and moment inequalities for U-statistics.}\\
 arXiv:math/0003228v1  [math.PR]  31 Mar 2000.\\

\vspace{3mm}

{\bf 9.  Ibragimov, R., Sharakhmetov, S.} {\it Exact bounds on the moments of symmetric statistics.}  Seventh Vilnius Conference on
Probability Theory and Mathematical Statistics, (1998). Abstracts of Communications, 243-244.\\

\vspace{3mm}

{\bf 10. Ibragimov, R., Sharakhmetov, S.} {\it On the exact Constants in the Rosenthal Inequality.}
Theory Probab. Appl., 1997, V. 42 pp. 294-302. \\

\vspace{3mm}

{\bf 11. Ibragimov, R., Sharakhmetov, S.} {\it The exact Constant in the Rosenthal Inequality for Sums of Independent real Random Variables with
Mean Zero.} Theory Probab. Appl., 2001, B.1 pp. 127-132. \\

\vspace{3mm}

{\bf 12. W. B. Johnson, G. Schechtman and J. Zinn.}  {\it Best constants in moment inequalities
for linear combinations of independent and exchangeable random variables. } Ann. Probab. 13 (1985), pp. 234-253.\\

\vspace{3mm}

{\bf 13. Johnson, W. B., Schechtman, G.} {\it Sums of independent random variables in rearrangement invariant function spaces.} Ann. Probab. 17 (1989), no. 2, 789-808.\\

\vspace{3mm}

{\bf 14. Johnson, W. B., Schechtman, G., Zinn, J.} {\it Best constants in moment inequalities for linear combinations of independent and exchangeable random variables.}
Ann.Probab. 13 (1985), no. 1, 234-253.\\

\vspace{3mm}

{\bf 15.  A.Khinchine.}  {\it Uber dyadische Br\"ucke.} Math. Z. 18 (1923), pp. 109-116. \\

\vspace{3mm}

{\bf 16.  R.Latala.}  {\it Estimation of moments of sums of independent real random variables. }
Ann. Probab. 25 (1997), pp. 1502-1513.\\

\vspace{3mm}

{\bf 17. Ostrovsky E. and Sirota L.} {\it Schl\"omilch and Bell series for Bessel's functions, with probabilistic applications.  }
arXiv:0804.0089v1  [math.CV]  1 Apr 2008 \\

\vspace{3mm}

{\bf 18.  V. H. de la Pena, R.Ibragimov and S.Sharakhmetov.} {\it  On extremal distribution and
sharp Lp-bounds for sums of multi-linear forms.}  Ann. Probab. 31 (2003), pp. 630-675.\\

\vspace{3mm}

{\bf 19. I.F.Pinelis and S. A.Utev.} {\it Estimates of the moments of sums of independent random variables.}
Theory Probab. Applications, {\bf 29,} (1984), 574-577.\\

\vspace{3mm}

{\bf 20. J.Pitman.} {\it Some probabilistic aspects of set partitions.} American. Math. Monthly. 104 (1997), pp. 201-209.\\

\vspace{3mm}

{\bf 21. J.Riordan.} {\it An Introduction to Combinatorial Analysis.} Wiley, New York, NY, 1980. \\

\vspace{3mm}

 {\bf 22. H.Robbins.} {\it A remark on Stirling's formula.}  American Math. Monthly 62 (1955), pp. 26-29.\\

\vspace{3mm}

{\bf 23. H.P.Rosenthal.} {\it On the subspaces of $ L(p) \ $ spanned by sequences of independent random variables.}
Israel J.Math., 8, (1970), 273-303.\\

\vspace{3mm}

{\bf 24. Gideon Schechtman.} {\it Extremal configurations for moments of sums of independent positive random variables.}
 Internet publication, PDF, February 12, 2007.\\

\vspace{3mm}

 {\bf 25. S.A.Utev.} {\it Extremal problems in moment inequalities.} Theory Probab. Applications,
{\bf 28}, (1984), 421-422. \\

\vspace{3mm}

{\bf 26. S.A.Utev.}  {\it Limit Theorems in Probability Theory.}  Transactions of Inst. Math., Novosibirsk, 1985,
pp. 56-75 (in Russian).\\

\vspace{3mm}

\end{document}